\documentclass[dvipsnames,11pt,a4paper,reqno]{amsart}

\usepackage[margin = 3cm]{geometry}

\usepackage[full]{textcomp}

\usepackage[utf8]{inputenc}

\usepackage{amstext, amsbsy, amssymb, amsmath, amsthm, amsfonts}

\usepackage[english]{babel}

\usepackage{bbm}

\usepackage{bm}

\usepackage{enumitem}

\usepackage{graphicx}

\usepackage{hyperref} 
\hypersetup{colorlinks=true, urlcolor= black, linkcolor=blue, citecolor=teal}

\usepackage{mathrsfs}

\usepackage{mathtools}

\usepackage{nicefrac}

\usepackage{subcaption}

\usepackage{import}
\usepackage{xifthen}
\usepackage{pdfpages}
\usepackage{transparent}

\usepackage{import}

\theoremstyle{plain}
\newtheorem{thm}{Theorem}
\newtheorem{lem}[thm]{Lemma}

\newtheorem{prop}[thm]{Proposition}

\theoremstyle{definition}
\newtheorem{definition}[thm]{Definition}

\newtheorem*{exmp*}{Example}

\theoremstyle{remark}
\newtheorem{rem}[thm]{Remark}
\newtheorem*{rem*}{Remark}

\mathtoolsset{showonlyrefs, showmanualtags}

\newcommand{\calC}{\mathcal{C}}
\newcommand{\calD}{\mathcal{D}}
\newcommand{\calE}{\mathcal{E}}
\newcommand{\calF}{\mathcal{F}}
\newcommand{\calG}{\mathcal{G}}

\newcommand{\calQ}{\mathcal{Q}}

\newcommand{\bbN}{\mathbb{N}}

\newcommand{\bbP}{\mathbb{P}}

\newcommand{\bbR}{\mathbb{R}}

\newcommand{\bbZ}{\mathbb{Z}}

\newcommand{\bfT}{\mathbf{T}}

\newcommand{\sfL}{\mathsf{L}}

\newcommand{\sfR}{\mathsf{R}}

\newcommand{\dt}{\partial}

\newcommand{\La}{\Lambda}
\newcommand{\om}{\omega}

\newcommand{\Ann}{\mathrm{Ann}}
\newcommand{\Rect}{\mathrm{Rect}}
\newcommand{\Exp}{\mathsf{Exp}}
\newcommand\Fin{\calF_{\mathrm{in}}}
\newcommand\Fout{\calF_{\mathrm{out}}}

\title{Generic arm separation for planar percolation models}
	
\author{Lo\"{i}c Gassmann}
\address{University of Fribourg, Switzerland}
\email{loic.gassmann@unifr.ch}

\author{Ioan Manolescu}
\address{University of Fribourg, Switzerland}
\email{ioan.manolescu@unifr.ch}

\author{Maran Mohanarangan}
\address{University of Fribourg, Switzerland and University of Innsbruck, Austria}
\email{maran.mohanarangan@unifr.ch}
	
\date{\today}

\begin{document}

\begin{abstract}
	Arm separation is an important technical tool in the study of critical and near-critical planar percolation systems. When combined with RSW estimates, it proves quasi-multiplicativity of arm-event probabilities, links the derivatives of different events and allows the identification of universal arm exponents. 
    
	We provide a very generic statement of arm separation that applies to arbitrary planar percolation models. Our argument is based on exploring configurations conditioned on the existence of arms in such a way that ``exposed'' arms are revealed with positive probability. The exploration algorithm is independent of the percolation model and applies to all polychromatic arm events. 
	As an application, we prove a strong arm-separation result for critical planar FK-percolation.
\end{abstract}

\maketitle

\section{Introduction}

Arm separation is a fundamental tool in the study of planar critical and near-critical percolation. An arm event denotes the existence of disjoint paths of prescribed colours joining the two boundaries of an annulus. Informally, arm separation states that, conditional on such an event, the arms can be chosen so that their endpoints on the outer  boundary are separated by distances comparable to the outer radius, with uniformly positive probability. An analogous statement applies at the inner boundary.

Arm separation was introduced by Kesten \cite{Kes87a} for planar Bernoulli percolation; see \cite{Nol08} for a detailed exposition. Its usefulness lies in the room it provides to extend arms and join them to other paths using Russo--Seymour--Welsh (RSW) estimates. Keeping the endpoints sufficiently far apart allows these connections to be constructed with uniformly positive probability, which is crucial when comparing events across different scales. Arm separation is thus a technical but essential ingredient in a wide range of applications, including quasi-multiplicativity \cite{Kes87a, Nol08, CheDumHon13, DumManTas21}, the derivation of scaling relations \cite{Kes87b,Kes87a,Nol08,GarPetSch13,GarPetSch18,DumMan20}, the construction of incipient infinite clusters \cite{Kes86, Jar03, DamSap11,BasSap17}, noise sensitivity \cite{SchSte10, GarPetSch10,TasVan25} and identification and comparison of arm exponents \cite{ SmiWer01,BefNol11,KohLeh25,GasMan25,RadTas26} to name just a few. 

In this paper, we define a related geometric property, which we call \emph{exposition}. An exploration algorithm will reveal the configuration step by step and, conditional on the arm event, will produce exposed paths of the prescribed colours with uniformly positive probability. These paths will have adjacent unexplored regions through which suitable crossings may be glued to them. The algorithm is the same for every percolation measure and requires no assumptions such as RSW estimates, positive association or independence.

Unlike classical proofs, which intertwine exploration and RSW techniques, our approach completely separates the geometric and probabilistic aspects of arm separation. We show how to explore the configuration in order to produce exposed arms. One may then use RSW estimates to glue these arms to desired destinations. 

As an illustration, we use our approach to produce a \emph{strong arm-separation result} for critical FK-percolation. 
This separation property applies to any polychromatic arm event (even non-alternating) in an annulus $\Lambda_R \setminus \La_{r}$ (with $r < R/2$) and will separate the arms on the outer boundary. It holds even conditionally on the configuration in $\La_{R/2}$. 
A precise statement is provided in Section~\ref{sec:applications_FK}.

\subsection{Definitions and main result}

A tiling of the plane is a family $\calG$ of closed, bounded sets, called \emph{cells}, with pairwise disjoint interiors and whose union is $\bbR^2$. 
We assume that every cell is simply connected and bounded by a simple closed curve of finite length.
Furthermore, we restrict our attention to locally finite tilings, meaning that every compact set of the plane intersects only a finite number of cells. 
Finally, we require that no point belongs to the boundaries of more than three cells. 

A percolation configuration on $\calG$ is an element of $\{0,1\}^\calG$, viewed as an assignment of colours to the cells of $\calG$. Cells of colour $0$ are called {\em closed}, while those of colour $1$ are called {\em open}; 
these will be represented by yellow and blue\footnote{The use of yellow and blue is common in the percolation community since the works of Smirnov~\cite{Smi01, CamNew06, BefDum13, KhrSkoSmi25}, who credited their choice to the flag of Sweden, where he was working at the time. We chose them with another country's flag in mind.}, respectively, in all figures. 
A percolation model $\bbP$ is a probability measure on the set $\{0,1\}^\calG$ equipped with the standard cylinder $\sigma$-algebra.
No further assumptions are imposed on $\calG$ or $\bbP$ at this stage.

Fix $0 < r < R/2$ and $z \in \bbR^2$.
Let $\La_r(z) = z+[-r,r]^2$ denote the box of radius~$r$ centred at $z$, and set $\La_r = \La_r(0)$. We write $\partial \La_r$ for its topological boundary and $\Ann(r,R) =\overline{ \La_R \setminus \La_{r}}$ for the closed annulus of inner radius $r$ and outer radius $R$. 
Finally, for any subset~$A \subset \bbR^2$, we write $\calG(A)$ for the union of cells of $\calG$ intersecting $A$. 

A colour sequence $\sigma$ is an element of $\{0,1\}^N$ for some $N \geq 1$. We call it \emph{polychromatic} if it contains both values $0$ and $1$, and \emph{monochromatic} otherwise. 

A cell-path $\gamma$ is a finite sequence of distinct cells in which consecutive cells are adjacent. Two cell-paths are disjoint if they have no cells in common. For a configuration $\om \in \{0,1\}^\calG$, a cell-path is called open (resp.\ closed) in $\om$ if all its cells are open (resp.\ closed) in $\om$.

\begin{figure}
\begin{center}
\includegraphics[width = 0.45\textwidth, page = 1]{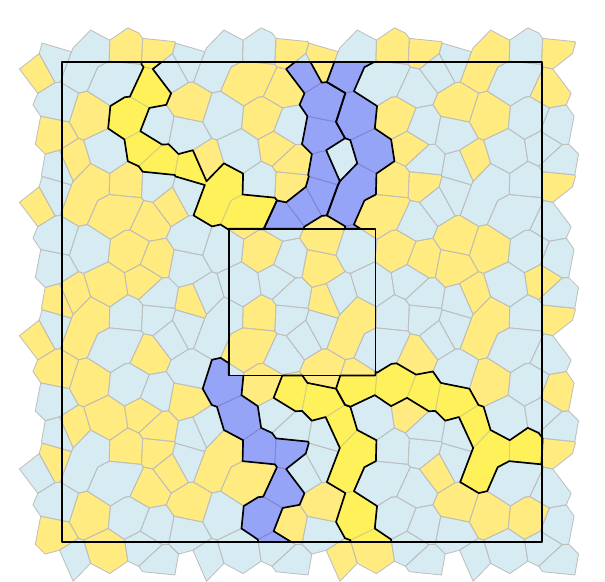}\hspace{ 0.05\textwidth}
\includegraphics[width = 0.45\textwidth, page = 2]{Tessellation.pdf}
\caption{{\em Left:} A percolation configuration on a tiling $\calG$. The annulus $\Ann(r,R)$ is depicted in black. 
The configuration contains six arms corresponding to the colour sequence $\sigma = (1,0,0,1,1,0)$. 
{\em Right:}
An exploration that exposes each arm. The exposing tubes are marked in red; the exposition zones are the unexplored regions surrounded by bold red lines.}
\label{fig:tessellation}
\end{center}
\end{figure}

Furthermore, a cell-path $\gamma$ in $\calG(\Ann(r,R))$ is said to be a crossing of the annulus $\Ann(r,R)$ if it contains a cell intersecting $\partial \La_r$ and a cell intersecting $\partial \La_R$. 
We say that $\Ann(r,R)$ contains arms of colours $\sigma$ if there exist disjoint cell-paths $\gamma_1, \dots, \gamma_N$ crossing $\Ann(r,R)$, indexed in counterclockwise order, such that $\gamma_i$ is open (resp.\ closed) if $\sigma_i\nobreak=\nobreak1$ (resp.\ $\sigma_i = 0$).
The paths $\gamma_i$ are called \emph{arms} and we write $A_{\sigma}(r,R)$ for the event corresponding to their existence. 
See Figure~\ref{fig:tessellation} (left) for an example. 

Our approach crucially relies on \emph{explored sets}. Loosely speaking, an explored set $\Exp$ is a union of cells of $\calG$ that is measurable in terms of the colours of the cells it contains, along with a source of randomness independent of the percolation configuration\footnote{We may need to enlarge the probability space on which the model is defined to accommodate general explored sets. We will do so without further mention.}. 
We will actually use the term ``explored'' in a more restrictive sense --- see Definition~\ref{def:Exp_set} below.

A tube of quality $\eta>0$ at scale $R$ is the union $T = \bigcup_{i = 1}^m \La_{\eta R}(x_i)$ where $x_1,\dots, x_m$ are distinct points belonging to a common translate of $2\eta R \bbZ^2$ with $\|x_i - x_{i-1}\|_2 = 2\eta R$ for all $2 \leq i \leq m$ and with $m \leq \eta^{-2}$. 
The starting box of $T$ is~$\La_{\eta R}(x_1)$ and its ending box is $\La_{\eta R}(x_m)$. 
A crossing of $T$ is any continuous curve in $T$ starting in $\La_{\eta R}(x_1)$ and ending in~$\La_{\eta R}(x_m)$.

Given an explored set $\Exp \subset \calG(\Ann(r,R))$, a cell-path $\gamma$ is said to be \emph{$\eta$-exposed towards the outside} if 
\begin{enumerate}[label = (\roman*)]
	\item $\gamma \subset \Exp$,
    \item there exists a tube $T$ of quality $\eta$ at scale $R$ contained in $\Ann(R/2,(1+2\eta)R)$, whose starting box lies in $\La_R^c$ and
    such that
    
	\item any crossing of $T$ necessarily intersects $\gamma$ and
	\item the connected component of the start of $T$ in $T \setminus \gamma$ --- henceforth referred to as the \emph{exposition zone} of $\gamma$ in $T$ --- is disjoint from $\Exp$.
\end{enumerate}
We will also say that $\gamma$ is exposed by $T$ in $\Exp^c$.
One can define a similar notion of \emph{$\eta$-exposition towards the inside} by considering tubes of quality $\eta$ at scale $r$, contained in $\Ann((1-2\eta)r,2r)$ with their starting box contained in $\La_r$. 

Finally, we say $\Exp$ \emph{$\eta$-exposes paths of colours $\sigma$ towards the outside} if
\begin{enumerate}[label = (\roman*)]
	\item $\Exp \subset \calG(\Ann(r,R))$,
	\item $\Exp$ contains disjoint cell-paths $\gamma_1,\dots, \gamma_N$ starting on $\partial \La_r$, with each $\gamma_i$ of colour $\sigma_i$, 
	indexed in counterclockwise order and
	\item the paths $\gamma_1,\dots, \gamma_N$ are $\eta$-exposed towards the outside, with pairwise disjoint exposition zones.
\end{enumerate}
See Figure~\ref{fig:tessellation} (right) for an illustration. 
A similar definition applies towards the inside. 

\begin{thm}[Arm exposition]\label{thm:arm_exposition}
	For any  polychromatic colour sequence $\sigma \in \{0,1\}^N$, there exist constants $\eta = \eta(\sigma)> 0$ and $c = c(\sigma) > 0$ such that the following holds.
	For any $0 < r < R/2$, any tiling of the plane $\calG$ such that all cells intersecting $\La_R$ have diameter at most $\eta R$,
	and any percolation measure $\bbP$ on $\calG$, there exists an explored set $\Exp$ such that 
	\begin{align}\label{eq:arm_exposition}
		\bbP\big[\text{$\Exp$ $\eta$-exposes paths of colours $\sigma$ towards the outside}\,\big|\, A_{\sigma}(r,R)\big] \geq c,
	\end{align}
	as long as the event in the conditioning has positive probability. An analogous result applies for exposition towards the inside.
\end{thm}

\begin{rem}
	A similar statement for monochromatic arm events is only available in a simply connected setting --- see Theorem~\ref{thm:arm_exposition_rect} below. 
	We were unable to prove such a statement for the annulus. However, in many applications involving monochromatic arm events, separation may be replaced by applications of the FKG inequality. 
\end{rem}

\subsection{Comments and applications}\label{sec:comments}

\subsubsection*{Minimal assumptions} It may be surprising that Theorem~\ref{thm:arm_exposition} applies to any percolation measure, with no assumption of properties of criticality, most notably the RSW property.
Note however that the exposition of the different arms is generally only useful when combined with versions of the RSW, Spatial Markov and FKG properties.

\subsubsection*{Uniformity over percolation models}
It will be apparent from the proofs that the set $\Exp$ is revealed by an exploration algorithm (formalised as a decision tree in Section~\ref{sec:exploration}) that queries the colours of cells one by one. When successful, the algorithm produces an explored set that $\eta$-exposes the desired paths. The same algorithm is used for all percolation models and its probability of success is bounded below by a positive constant depending only on the colour sequence $\sigma$. This allows us to state the theorem uniformly over percolation models on tilings. 

\subsubsection*{Probability of success}
The probability of success in \eqref{eq:arm_exposition} may be rendered arbitrarily close to 1 by shrinking $\eta > 0$. Indeed, for fixed $\delta > 0$, taking $\eta > 0$ small enough, one may force the exploration in Theorem~\ref{thm:arm_exposition} to satisfy
\begin{align}
    \Exp \subset \calG(\Ann(r,(\alpha+\delta) R))
\end{align}
for any $1/2 \leq \alpha \leq 1-\delta$, without affecting the probability of success $c$. Repeating the procedure in successive annuli of ``thickness'' $\delta R$  yields independent attempts at exposition, so that the overall probability of success is at least $1-(1-c)^{1/(2\delta)}$. Taking  $\delta$ small enough renders the probability of success arbitrarily close to one.

\subsubsection*{Robustness}
\begin{figure}
	\begin{center}
	\includegraphics[width = 0.17\textwidth, page = 1]{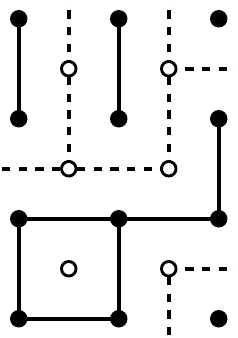}\hspace{ 0.03\textwidth}
	\includegraphics[width = 0.17\textwidth, page = 2]{cells_square.pdf}\hspace{ 0.15\textwidth}
	\includegraphics[width = 0.17\textwidth, page = 1]{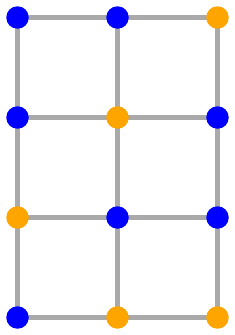}\hspace{ 0.03\textwidth}
	\includegraphics[width = 0.17\textwidth, page = 2]{cells_square2.pdf}
	
	\caption{{\em Left:} A bond percolation on the square lattice and its interpretation as a colouring of cells. 
	The cells corresponding to vertices are all blue (open), while those corresponding to faces are all yellow (closed). The cells corresponding to edges are coloured depending on the state of the corresponding edge. 
	{\em Right:} A site percolation on the square lattice and its interpretation on the same graph of cells $\calG$. Here the cells corresponding to edges are all blue and those corresponding to faces are all yellow. The cells corresponding to vertices carry all the randomness.}
	\label{fig:lattices_as_tilings}
	\end{center}
\end{figure}

The interpretation of percolation as colourings of cells allows one to incorporate site or bond percolation models on any bi-periodic planar lattice, and much more. Bond percolation is recovered by viewing the colour $0$ as the dual percolation model, while site percolation is recovered by viewing closed cells as the closed vertices, with the connectivity notion given by the matching graph.
See Figure~\ref{fig:lattices_as_tilings} for examples of cell representations of bond and site percolation on $\bbZ^2$.

Theorem~\ref{thm:arm_exposition} also applies to less classical settings, such as Voronoi percolation or percolation on random planar maps, and, with a natural adaptation, to Boolean percolation, where balls of random radius are placed in the plane according to a Poisson point process. In the latter setting, arm events for alternating colour sequences are defined exactly as above. A separate definition is needed to make sense of several consecutive arms in the occupied set, while the corresponding notion for the vacant set does not make sense.

We believe that those interested in applying Theorem~\ref{thm:arm_exposition} in these non-standard settings will find it easily adaptable.

One caveat applies to models such as Voronoi or Boolean percolation, where the tiling itself is random. Theorem~\ref{thm:arm_exposition} may be generalised to this case with $\Exp$ revealing both the colour \emph{and} the shape of each explored cell. We will not detail this here, as the adaptation is immediate. We do however warn the interested reader that, in such contexts, the knowledge of the cell shapes in $\Exp$ may affect the RSW estimates in the unexplored regions. For more background on RSW estimates in Voronoi and Boolean percolation, see \cite{Tas16, AhlGriMorTas16} and \cite{AhlTasTei18}, respectively. It is worth mentioning that the quenched result of \cite{AhlGriMorTas16} may actually allow one to perform arm exposition by first revealing the Voronoi tessellation and then applying the algorithm of Theorem \ref{thm:arm_exposition} on it.

\subsubsection*{Motivation} 
The present work is mainly motivated by applications to critical FK-percolation, where RSW techniques give lower bounds on the probability of open and closed crossings of tubes of given quality. Combined with the FKG inequality and the spatial Markov property, these estimates allow one to construct crossings in unexplored regions and glue them to exposed paths.
This is the basic mechanism underlying arm separation.

The classical proofs of arm separation for Bernoulli percolation \cite{Kes87a,Nol08} extend to FK-percolation when suitable RSW estimates are available. For alternating arm events, the usual box-crossing estimate (see \eqref{eq:RSW}) is sufficient. For general colour sequences, a complete proof for FK-Ising percolation is given in \cite{CheDumHon13}, using stronger RSW estimates for arbitrary topological rectangles. The same approach extends to $1 \leq q < 4$ using \cite{DumManTas21}, but does not cover $q = 4$, where these stronger estimates are expected to fail.

Our argument treats arbitrary polychromatic colour sequences using only RSW estimates uniform over boundary conditions at macroscopic distance. These estimates remain valid at $q = 4$, so our approach provides a unifying framework.
Moreover, it yields a stronger arm-separation result which we describe next.

\subsubsection*{Strong arm separation under conditioning}
An important consequence of our approach is that arm separation remains valid under arbitrary compatible conditioning in the inner annulus $\Ann(r,R/2)$. Indeed, fix a configuration in $\Ann(r,R/2)$ under which $A_{\sigma}(r,R)$ has positive conditional probability. Since Theorem~\ref{thm:arm_exposition} applies uniformly over percolation measures, it also applies to the FK-percolation measure conditioned on the fixed configuration in $\Ann(r,R/2)$.

To turn the exposed paths into a more standard notion of well-separated arms (see Definition~\ref{def:separation}), we use RSW constructions in the exposition zones. As we will show in the proof, these zones lie at macroscopic distance from the conditioned region. The spatial Markov property and the uniformity of RSW estimates over boundary conditions at such distances therefore give a positive lower bound on the probability of these constructions, independently of the prescribed configuration in $\Ann(r,R/2)$. This yields arm separation uniformly under this conditioning; see Theorem~\ref{thm:arm_separation} for the precise statement.

In particular, our approach implies the ``super-strong'' separation statement conjectured by Garban, Pete and Schramm in \cite[Conjecture~A.4]{GarPetSch13}. A half-plane version was previously obtained in \cite{DuGao+24} for critical site percolation on the triangular lattice. 
Our result covers the full-plane case and extends this form of separation to critical FK-percolation.

\subsection{Exposed crossings in rectangles}

An intermediate step in our proof is an exposition result similar to Theorem~\ref{thm:arm_exposition}, but for crossings of quads (that is, topological rectangles). We state it in the context of Euclidean rectangles for simplicity. 

For $H \geq 2R > 0$, define $\Rect(R, H) = [-R,R] \times [0,H]$, with its four sides denoted $\partial_T \Rect(R, H)$, $\partial_L \Rect(R, H)$, $\partial_B \Rect(R, H)$ and $\partial_R \Rect(R, H)$, in counterclockwise order, starting from the top. 
For a colour sequence $\sigma$, let $ A_\sigma(\Rect(R,H))$ be the event that~$\Rect(R,H)$ contains disjoint cell-paths 
from $\partial_{\mathrm B}\Rect$ to $\partial_{\mathrm T}\Rect$ of colours $\sigma_1, \dots, \sigma_N$ when indexed from left to right. 

An explored set $\Exp \subset \calG(\Rect(R,H))$ is said to \emph{$\eta$-expose a cell-path $\gamma$ towards the top} if 

\begin{enumerate}[label = (\roman*)]
	\item $\gamma \subset \Exp$,
    \item there exists a tube $T$ of quality $\eta$ at scale $R$ whose starting box is contained in $[-R,R] \times [H, +\infty)$,
	\item any crossing of $T$ necessarily intersects $\gamma$ and
	\item the connected component of the start of $T$ in $T \setminus \gamma$ --- henceforth referred to as the exposition zone of $\gamma$ in $T$ --- is entirely contained in the unexplored part of $[-R,R] \times [H - R, +\infty)$. 
\end{enumerate}
As before, we will also say that $\gamma$ is exposed by $T$ in $\Exp^c$. 
That the exposition zones are required to be above the height $H-R$ is similar to them being required to be outside of~$\La_{R/2}$ in Theorem~\ref{thm:arm_exposition}. It ensures that the exposition actually occurs near the top boundary of the rectangle.

Finally, we say $\Exp$ \emph{$\eta$-exposes crossings of colours $\sigma$ towards the top} if
\begin{enumerate}[label = (\roman*)]
	\item $\Exp \subset \calG(\Rect(R,H))$,
	\item $\Exp$ contains disjoint cell-paths $\gamma_1,\dots, \gamma_N$ starting on $\partial_B \Rect(R,H)$ with each $\gamma_i$ of colour $\sigma_i$, 
	indexed from left to right and
	\item the cell-paths $\gamma_1,\dots, \gamma_N$ are $\eta$-exposed towards the top, with their exposition zones being disjoint. 
\end{enumerate}

\begin{thm}\label{thm:arm_exposition_rect}
	For any colour sequence $\sigma \in \{0,1\}^N$, there exist constants $\eta = \eta(\sigma)> 0$ and $c = c(\sigma) > 0$ such that the following holds.
	For any  $H \geq 2R > 0$, any tiling $\calG$ such that all cells intersecting $\Rect(R,H)$ have diameter at most $\eta R$ and 
	any percolation measure~$\bbP$ on $\calG$, there exists an explored set $\Exp$ such that
	\begin{align}\label{eq:arm_exposition_rect}
		\bbP\big[\text{$\Exp$ $\eta$-exposes crossings of colours $\sigma$ towards the top}\,\big|\, A_{\sigma}(\Rect(R,H))\big] \geq c,
	\end{align}
	as long as the event in the conditioning has positive probability. 
\end{thm}

Unlike Theorem~\ref{thm:arm_exposition}, this result also applies to monochromatic colour sequences.
Indeed, the fact that the rectangle has natural notions of left and right boundaries --- which boils down to it being simply connected --- renders the proof simpler. Actually, we deduce Theorem~\ref{thm:arm_exposition} by cutting the annulus so that it becomes simply connected without altering the arm event, then applying Theorem~\ref{thm:arm_exposition_rect}. In the cutting step, we use that the colour sequence is polychromatic to explore an interface.

\vspace{10pt}

\noindent
\textbf{Structure of the paper.} In Section~\ref{sec:exploration}, we formally define the notion of explored set. 
Section~\ref{sec:proof_box} contains the proof of Theorem~\ref{thm:arm_exposition_rect}, which is then used in Section~\ref{sec:proof_annulus} to show Theorem~\ref{thm:arm_exposition}. 
Finally, Section~\ref{sec:applications_FK} contains consequences for FK-percolation.

\section{Exploration algorithms}\label{sec:exploration}

As stated in the introduction, the explored set $\Exp$ will be the result of an exploration procedure. 
Such a procedure consists of revealing a random configuration $\omega$ cell by cell, with the choice of the next cell to be revealed depending on all previously revealed information and, potentially, on an external source of randomness that is independent of $\omega$. 
We formalise this notion using \emph{decision trees}.

\subsection{Decision trees}

Fix a tiling $\calG$ and a configuration $\omega$. For a sequence of cells $e = (e_1, e_2,\dots)$ and $t \geq 0$, we write $e_{[t]} = (e_1, \dots, e_t)$ and $\om_{e_{[t]}} = (\om_{e_1}, \dots, \om_{e_t})$. Write $\calG_{[t]}$ for the collection of ordered $t$-tuples of cells, so that $e_{[t]} \in \calG_{[t]}$. We will sometimes slightly abuse notation and also write $e_{[t]}$ for the set $\{e_1, \dots, e_t\}$ whenever it is convenient.

A decision tree $\bfT$ is an initial cell $e_1$ and a set of functions $\psi_t : \calG_{[t]} \times \{0,1\}^t \to \calG$ for~$t \geq 1$. 
It should be seen as associating to each pair $(e_{[t]}, \omega_{e_{[t]}})$ of revealed cells the next cell to reveal. 
In particular, we ask that $\psi_t (e_{[t]}, \omega_{e_{[t]}}) \notin \{e_1, \dots, e_t\}$. 

When the decision tree is fixed, we write $(e_{[t]}, \omega_{e_{[t]}})_{t \geq 1}$ for the unique sequence starting with $(e_1,\omega_{e_1})$ and such that 
$e_{t+1} = \psi_t (e_{[t]}, \omega_{e_{[t]}})$ for all $t \geq 1$. 
In the rest of the paper, we will interpret decision trees as exploration algorithms and refer to them as such.

A \emph{stopping time} for $\bfT$ is a random variable $\tau$ taking values in $\bbN \cup \{\infty\}$ which is such that $\{\tau \leq t\}$ is measurable with respect to $(e_{[t]},\om_{e_{[t]}})$. The set of cells explored up to the stopping time $\tau$ (including the order in which they were revealed) is $e_{[\tau]}$ and the configuration revealed at that time is $\omega_{e_{[\tau]}}$.

We call a decision tree \emph{random} if both the choice of $e_1$ and the functions $\psi_t$ also depend on an external source of randomness, which may be encoded by a sequence $(U_t)_{t \in \bbN}$ of independent uniform variables on $[0,1]$ which is independent of $\omega$.
The notion of stopping time adapts readily, by imposing that $\{\tau \leq t\}$ is measurable with respect to $(e_{[t]},\om_{e_{[t]}},(U_s)_{s \leq t})$.

\begin{definition}\label{def:Exp_set}
    An explored set $\Exp$ is a random set of cells of $\calG$ for which there exists a random decision tree and a stopping time $\tau$ such that $\Exp = e_{[\tau]}$.
\end{definition}

\subsection{Interface exploration}

\begin{figure}
\begin{center}
\includegraphics[width = 0.4\textwidth]{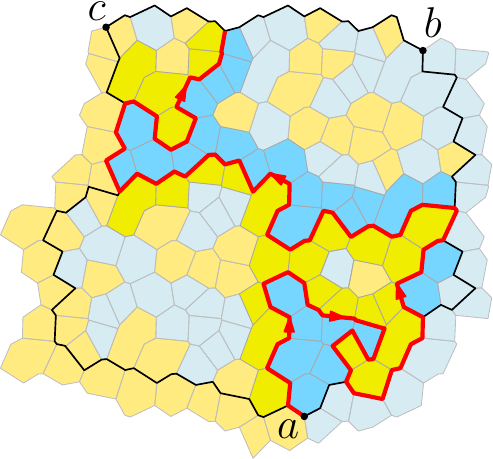}
\caption{\emph{Interface exploration:} A domain $\Omega$ surrounded by the black curve, with three points $a,b,c$ on the boundary. The exploration path starting at $a$ and ending on the arc $(b,c)$ is depicted in red. The cells of $\Omega^c$ adjacent to the arcs $(a,b)$ and $(c,a)$ are considered open and closed, respectively. Only the cells in $\Omega$ adjacent to the exploration path are explored (brighter colours).}
\label{fig:exploration}
\end{center}
\end{figure}

An important example is the exploration of an interface between open and closed cells. See Figure~\ref{fig:exploration} for a helpful illustration.

Let $\Omega$ be a finite union of cells whose boundary is a simple closed curve.
Let $a,b$ and $c$ be three distinct points on the boundary of $\Omega$, placed in counterclockwise order, each at the interface of three cells, with exactly one of the three cells belonging to $\Omega$. 
Write~$(a,b)$,~$(b,c)$ and $(c,a)$ for the three boundary arcs going in counterclockwise direction.

The exploration curve $\Gamma$ will be a path of edges of $\calG$ started from $a$ such that all cells on one side of $\Gamma$ are of one colour and all cells on the other side of $\Gamma$ are of the opposite colour. For definiteness, we describe an exploration following edges with closed on the left and open on the right. To ensure that $\Gamma$ remains inside $\Omega$, we impose open boundary conditions on the arc $(a,b)$ and closed boundary conditions on $(c,a)$. This will be made precise below. The remaining arc, including its endpoints $b$ and $c$, is the target at which the exploration stops.

We write $(\Gamma_n)_{n \geq 0}$ and $(\Gamma_{(n,n+1)})_{n \geq 0}$ for the sequence of vertices and oriented edges\footnote{Note that the index $n$ is distinct from the number of revealed cells $t$ in the exploration.}, respectively, visited by $\Gamma$. Set $\Gamma_0 = a$. To properly initialise the procedure, we also define~$\Gamma_{(-1,0)}$ as the unique edge between the two exterior cells incident to $a$, oriented towards~$a$.

At each step $n \geq 0$, first check whether $\Gamma_n$ lies on the target arc $(b,c)\cup\{b,c\}$. If so, stop the exploration. Otherwise, proceed as follows.

For $n = 0$, let $e$ be the unique cell in $\Omega$ incident to $a$. For $n \geq 1$, let $e$ be the unique cell incident to $\Gamma_n$ that is not incident to the incoming edge $\Gamma_{(n-1,n)}$. 
If $e$ belongs to~$\Omega$, reveal its colour (unless it has already been revealed) and let $\Gamma_{(n,n+1)}$ be the unique outgoing edge adjacent to $e$ that preserves closed on the left and open on the right. If $e$ lies outside~$\Omega$, use the boundary condition imposed above. 
More precisely, the following rules should be applied for $n \geq 0$:
\begin{enumerate}[label=(\roman*)]
    \item if $e \subset \Omega$ and $e$ is open, let $\Gamma_{(n,n+1)}$ be the outgoing edge to the left as seen from $\Gamma_{(n-1,n)}$;
    \item if $e \subset \Omega$ and $e$ is closed, let $\Gamma_{(n,n+1)}$ be the outgoing edge to the right as seen from $\Gamma_{(n-1,n)}$;
    \item if $e \not \subset \Omega$ and $\Gamma_n$ lies on $(a,b)$, pretend that $e$ is open and follow the same procedure as in (i);
    \item if $e \not \subset \Omega$ and $\Gamma_n$ lies on $(c,a)$, pretend that $e$ is closed and follow the same procedure as in (ii).
\end{enumerate}
After choosing the outgoing edge $\Gamma_{(n,n+1)}$, denote its endpoint by $\Gamma_{n+1}$. 

Recording the cells revealed during this procedure defines a decision-tree exploration. Due to our imposition of boundary conditions and $\Omega$ containing only finitely many cells, the exploration curve will reach the target after finitely many steps. In particular, the number $\tau$ of cells that are revealed before the exploration stops is finite and is a stopping time for the associated decision tree. Formally, the exploration is only defined up to $\tau$, but one may continue the exploration after $\tau$ using arbitrary rules.

\section{Exposition in a rectangle: Proof of Theorem~\ref{thm:arm_exposition_rect}}\label{sec:proof_box}

Before delving into the proof of Theorem~\ref{thm:arm_exposition_rect}, let us make some useful observations about the notion of tubes.

\subsection{Cutting tubes}\label{sec:cutting_tubes}
Recall that a tube $T$ of quality $\eta$ at scale $R$ is $T = \bigcup_{i=1}^m \Lambda_{\eta R}(x_i)$ for a sequence $x_1, \dots, x_m \in v + 2\eta R \bbZ^2$ of distinct points with $\Vert{x_i-x_{i-1}}\Vert_2 = 2\eta R$ for all~$i = 2, \dots m$ and $1 \leq m \leq \eta^{-2}$ and some fixed translation vector $v \in \bbR^2$.

\begin{figure}
    \begin{center}
    \includegraphics[width = .4\textwidth, page = 1]{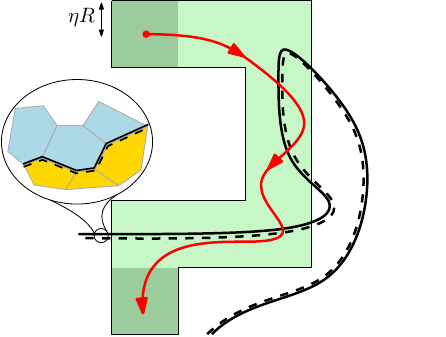}
    \hspace{ .05\textwidth}
    \includegraphics[width = .4\textwidth, page = 2]{cut_tube.pdf}
    \caption{
        {\em Left:} A tube $T$ of quality $\eta$ at scale $R$ (in green) with the start and end marked by shaded squares. 
        An exploration curve (in black) with the open side marked by a solid line and the closed side marked by a dashed line. 
        Its open side is exposed in $T$. The red crossing of the tube hits the open side of the exploration curve first. 
        {\em Right:} The lengthwise cutting of~$T$ into two tubes of quality $\eta / 2$. Each square $\La_{\eta R}(\cdot)$ is cut depending on the position of the previous and following box.
    }
    \label{fig:cut_tube} 
    \end{center}
\end{figure}

We may cut any tube of quality $\eta$ at scale $R$ lengthwise into two tubes $T_1, T_2$ of quality~$\eta/2$ such that $T_1 \cup T_2 = T$ and $T_1$ and $T_2$ have disjoint interiors. This is done as follows (see Figure~\ref{fig:cut_tube} (\emph{right}) for an example). 

Denote by $\zeta$ the polygonal path obtained by joining each consecutive pair of points $x_i$ and $x_{i+1}$ by a line segment for all $i = 1, \dots, m-1$. Extend the first and last segments to the boundaries of the corresponding starting and ending boxes such that they are cut in half.
Note that $\zeta$ naturally cuts each box $\Lambda_{\eta R}(x_i)$ of $T$ in two parts $P_1 (x_i)$ and $P_2 (x_i)$ that lie at its right and at its left, respectively. 
For each $i$, $P_1 (x_i)$ and $P_2 (x_i)$ are each a union of one, two or three boxes $\Lambda_{\eta R/2}(z)$ for $z \in  v + (\tfrac{\eta R}{2}, \tfrac{\eta R}{2})+\eta R \bbZ^2$. 
On each side of~$\zeta$, the subsquares can be ordered along the path to form a nearest-neighbour chain. Each chain contains at most $3m$ squares and $3m \leq (\eta/2)^{-2}$.
Hence, we obtain tubes of quality~$\eta/2$ by taking $T_1 = \bigcup P_1(x_i)$ and $T_2 = \bigcup P_2(x_i)$.

An immediate consequence of this construction is that if a cell-path $\gamma$ is exposed by~$T$ in $\Exp^c$, then $\gamma$ is also exposed by $T_1$ and $T_2$ in $\Exp^c$. Indeed, any crossing of $T_1$ is in particular a crossing of $T$, and hence intersects $\gamma$; and the connected component of the start of $T_1$ in $T_1 \setminus \gamma$ is contained in the connected component of the start of $T$ in $T \setminus \gamma$, which is unexplored by assumption. The same arguments apply to $T_2$.

\subsection{Proof of Theorem~\ref{thm:arm_exposition_rect} for favourable tilings}\label{sec:favourable_tilings}

Fix a tiling $\calG$ and a percolation measure $\bbP$. Also, fix $R,H$ and $\sigma$ as in the statement of the theorem. 
We start by proving the result in the case where all cells of $\calG$ have diameter at most $1$ and where $R$ is assumed  to be large enough.
Below, we will define different tubes of quality $\eta$ at scale $R$. Since $R$ and even $\eta R$ will be taken large, we may ignore cells crossing the boundary of the tubes. 
More precisely, we will suppose that $\Rect(R,H)$ and all tubes considered hereafter can be approximated by unions of cells.
We will discuss these approximations and the case of more general tilings $\calG$ at the end of the proof. 

A quad is a simply connected domain $\calQ$ with four marked points $a,b,c,d$ on its boundary in counterclockwise order. 
This defines four boundary arcs $(a,b)$, $(b,c)$, $(c,d)$ and $(d,a)$, which we also call the bottom, right, top and left of $\partial \calQ$, respectively. 
Below, all quads will be unions of cells, so one can perform the interface exploration described in the previous section in these domains.

\begin{proof}[Proof of Theorem~\ref{thm:arm_exposition_rect}]
	Recall that $N$ is the length of the colour sequence $\sigma$. 
	Since the scale~$R$ at which we work is fixed, we will omit it from notation. In particular we henceforth 
	write $\eta$-exposed and tube of quality $\eta$ without referring to the scale $R$. 
	
	Since the result is stated with no assumptions on the percolation measure, we may work under $\bbP[\, \cdot \,|\, A_{\sigma}({\rm Rect}(R,H))]$.
	
	We will construct $\Exp$ by progressively revealing interfaces up to different stopping times $\tau_1 \leq \dots \leq \tau_N$. 
	At each stage, we may have $\tau_k = \infty$, which means that the exploration procedure failed. 
	The procedure is successful if $\tau_N < \infty$ and we set $\Exp = e_{[\tau_N]}$. 
	
    Set $\eta_{k} = 2^{-k-3}$ for $0\leq k \leq N$. 
	At the end of stage $k \leq N$, if $\tau_k < \infty$ then the following conditions are satisfied.
    
	\begin{enumerate}
		\item There exist indices $\ell$ and $r$ with $\ell+N+1-r = k$ such that the explored region contains cell-paths $\gamma_1, \dots, \gamma_\ell$ and $\gamma_{r}, \dots, \gamma_N$ of colours $\sigma_1, \dots, \sigma_\ell$ and $\sigma_{r}, \dots, \sigma_N$, respectively, that are $\eta_k$-exposed. 
		Write $T_1^{(k)},\dots, T_\ell^{(k)}, T_{r}^{(k)}, \dots, T_{N}^{(k)}$ for the tubes of quality $\eta_k$ exposing the paths $\gamma_i$; they will all be defined on a common translated grid.	
		
		If $\ell = 0$, we will also have a tube $T_0^{(k)}$ that is disjoint from the above and which exposes\footnote{\label{footnote_curve_expo}The notion of exposition was only defined for cell-paths and not for general paths in $\bbR^2$. In this case, we simply mean that conditions (ii), (iii) and (iv) of the definition of exposition are satisfied for the corresponding segment.} the segment $\{-R\} \times[0,H - 4\eta_k R]$. 
		Similarly, when $r = N+1$, we have an additional tube $T_{N+1}^{(k)}$, disjoint from the above, and which exposes the segment $\{R\} \times[0,H - 4\eta_k R]$. 
		
		\item The explored part at this stage satisfies
		$e_{[\tau_k]} \subset	 {\rm Rect}(R,H - 4\eta_k R)$ and 
		the unexplored part $e_{[\tau_k]}^c := {\rm Rect}(R,H) \setminus e_{[\tau_k]}$ is a simply connected domain.
	
		Moreover, when $k < N$, the intersection of $\partial e_{[\tau_k]}^c$ with the bottom of $\Rect(R,H)$ is exactly the segment of the latter contained between the starting cells of $\gamma_\ell$ and~$\gamma_{r}$. We call this the bottom of $e_{[\tau_k]}^c$.
		
		\item There exist a.s.\ disjoint cell-paths from the bottom of $e_{[\tau_k]}^c$ to the top of $\Rect(R,H)$, avoiding $\gamma_\ell$ and $\gamma_{r}$, ordered from left to right, of colours $\sigma_{\ell + 1},\dots, \sigma_{r-1}$.
		These paths may use both unexplored and explored cells, but no cells of $\gamma_\ell$ and $\gamma_{r}$.
	\end{enumerate}
	See Figures~\ref{fig:T_0} and \ref{fig:stage_k} for illustrations.
    
	\begin{figure}
    	\begin{center}
    	\includegraphics[width = 0.6\textwidth]{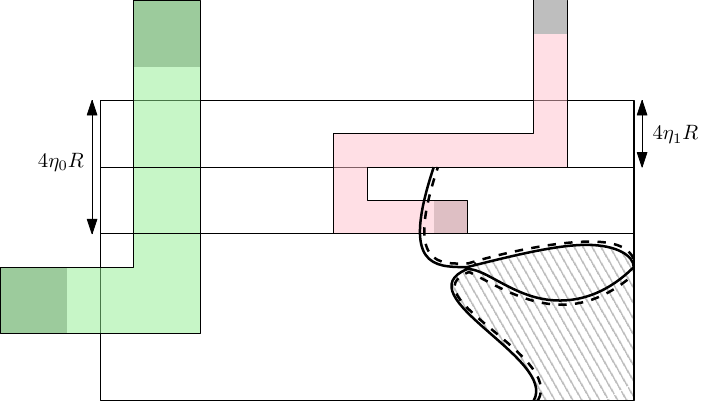}
    	\caption{The green tube is $T_0^{(0)}$; its start and end are shaded. This tube exposes the left side of the rectangle.
    	On the right, a primal path has been explored. The explored region at this stage is hatched; 
    	the red tube is $T_N^{(1)}$}
    	\label{fig:T_0}
    	\end{center}
	\end{figure}
	
	We will construct the exploration inductively and simultaneously check the points (1), (2) and (3) above. 
	At the end, we will bound the probability that $\tau_N <\infty$. 	
	As mentioned before, an external source of randomness will be used during the procedure. 
	This will be encoded by a series of i.i.d.\ uniform variables $X_1,\dots, X_N$ in $\{\sfL,\sfR\}$. 
	
	\medskip 
	
	\noindent {\bf Initial step:} To start, set $\tau_0 = 0$, $e_{[\tau_0]} = \varnothing$ and let
	\begin{align}
		T_0^{(0)} =& [-(1- \eta_0) R, -(1 - 3 \eta_0)R]\times [H - 7 \eta_0R,H+ 3\eta_0 R]  \\&\cup [-(1 + 3\eta_0)R,-(1- \eta_0) R] \times [H - 7 \eta_0R,H - 5\eta_0 R]
		\qquad \text{ and } \\
		T_{N+1}^{(0)} =& [(1-3 \eta_0) R, (1 - \eta_0)R]\times [H - 7 \eta_0R,H+ 3\eta_0 R] \\&  \cup [(1 -\eta_0)R,(1+3 \eta_0) R] \times [H - 7 \eta_0R,H - 5\eta_0 R],
	\end{align} 
	be $L$-shaped tubes; see Figure \ref{fig:T_0} for an example. These two tubes expose the left and right side of the rectangle, respectively. They only serve to initialise the induction process.
	Set $\ell = 0$ and $r = N+1$.
	All points are clearly satisfied when $k = 0$. 
    For completeness, we set $\gamma_0 = \{-R\} \times[0,H - 4\eta_k R]$ and $\gamma_{N+1} = \{R\} \times[0,H - 4\eta_k R]$. 
    This is a slight abuse of notation since the other $\gamma_i$ will refer to cell-paths.
	
	\medskip 
	
	\noindent {\bf Induction step: construction.} Let $0 \leq k < N$ and consider the exploration up to stage $k$, assuming that $\tau_{k} < \infty$. In particular, this means that conditions (1), (2) and (3) are satisfied at step $k$. 
    Proceed as follows for step $k+1$ of the algorithm.  
    
	Cut the tubes $T_{\ell}^{(k)}$ and $T_{r}^{(k)}$ lengthwise. 
	Write $T_{\ell}^{(k+1)}$ for the left half of the former and $S_L$ for its right half. 
	Similarly, write $T_{r}^{(k+1)}$ for the right half of the latter and $S_R$ for its left half. 
	For $j < \ell$ and $j > r$, cut the tube $T_{j}^{(k)}$ lengthwise and set $T_j^{(k+1)}$ to be the left half if $j < \ell$ and the right half if $j > r$\footnote{This choice is arbitrary. Indeed, the other halves would work just as well by the discussion in Section~\ref{sec:cutting_tubes}.}. 
    Observe that the tubes $T_{j}^{(k+1)}$ for $j \leq \ell$ and $j \geq r$ are indeed of quality $\eta_{k+1} = \eta_k/2$.
    
	\begin{figure}
    	\begin{center}
    	\includegraphics[width = 0.75\textwidth, page = 1]{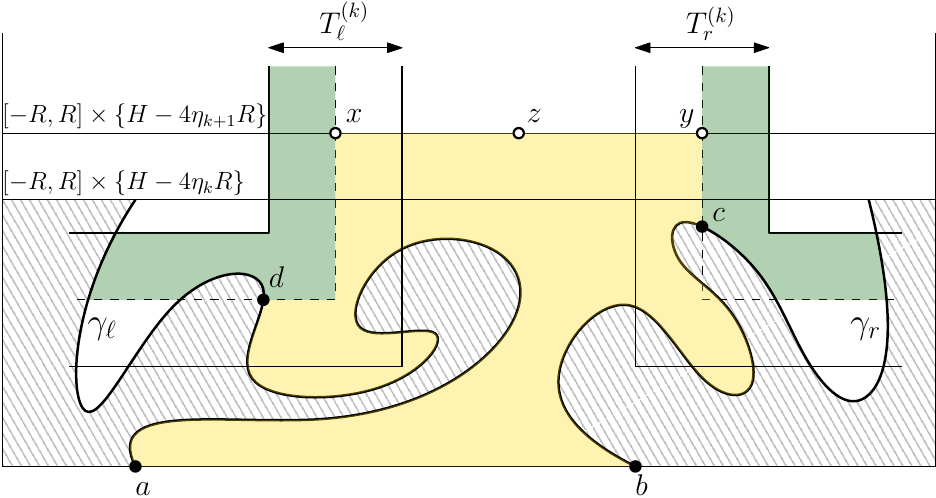}
    	\caption{An example of an exploration (hatched) at stage $k$ assuming $\tau_k < \infty$. The black paths are $\gamma_\ell$ and $\gamma_r$; 
    	they are exposed through the tubes $T_\ell^{(k)}$ and $T_r^{(k)}$, respectively. 
    	These tubes are cut lengthwise into $T_\ell^{(k+1)}$ and $T_r^{(k+1)}$ (green) and $S_L$ and $S_R$. 
    	The quad $\calQ_{k+1}$ is depicted in yellow.}
    	\label{fig:stage_k}
    	\end{center}
	\end{figure}

	We now define a quad $\calQ_{k+1}$; see Figure~\ref{fig:stage_k} for an illustration. 
	Consider the curve $\chi_L$ obtained by following the left boundary of $S_L$ from the top of $\Rect(R,H)$ up to the first time it reaches the boundary of a cell of $\gamma_\ell$ --- denote this point by $d$. 
	Then, continue extending $\chi_L$ by following the boundary of $\gamma_\ell$ down to the first time it reaches the bottom of $\Rect(R,H)$. Observe that the curve $\chi_L$ connects top and bottom of $\Rect(R,H)$. We write $a$ for its endpoint on the bottom of $\Rect(R,H)$.
	
	Let $\chi_R$ be defined similarly, using the right boundary of $S_R$ and $\gamma_r$ instead of $\gamma_\ell$. Denote by $c$ the point where the boundary of $S_R$ meets the boundary of a cell of $\gamma_r$ for the first time and by $b$ the endpoint of $\chi_R$ on the bottom of ${\rm Rect}(R,H)$.
	The fact that the boundaries of the tubes intersect $\gamma_\ell$ and $\gamma_r$, respectively, is guaranteed by their respective exposition. 
	 
	By point (2), the horizontal segment $(a,b)$ is exactly the bottom part of $\partial e_{[\tau_k]}^c$.
	Let~$x$ be the last point of $\chi_L$ on $[-R,R] \times \{H - 4 \eta_{k+1} R\}$ when $\chi_L$ is considered as a curve oriented towards $a$.
	Similarly, let $y$ be the last intersection of $\chi_R$ with this horizontal line. 
	Write~$(x,y)$ for the horizontal segment between $x$ and $y$. 
	Let $\calQ_{k+1}$ be the region delimited by the horizontal segment $(a,b)$, the sub-arc $(b,y)$ of $\chi_R$, the horizontal segment~$(y,x)$
	and the sub-arc $(x,a)$ of $\chi_L$; see Figure~\ref{fig:stage_k}.
	By closely following the definitions of $\chi_L$ and~$\chi_R$, one may notice that the cells of $\gamma_\ell$ and $\gamma_r$ are outside of $\calQ_{k+1}$; some of them being adjacent to the arcs $(d,a)$ and $(b,c)$ of the quad $\partial \calQ_{k+1}$. 
	
	Note that $(x,y)$ has length at least $2\eta_{k+1} R = \eta_k R$. Write $z$ for its midpoint. 
	The points $a,b,c,d$ defined above lie  on the boundary of $\calQ_{k+1}$ and are distributed in counterclockwise order. 
	The arc $(c,d)$ of $\partial \calQ_{k+1}$ contains the segment $(y,x)$. 
	
	As a consequence of point (3), 
	the quad $\calQ_{k+1}$ is crossed from $(a,b)$ to $(c,d)$ by disjoint paths of colours $\sigma_{\ell+1},\dots, \sigma_{r-1}$ ordered from left to right.

	If $X_{k+1} = \sfL$, we set $\gamma_{\ell+ 1}$ to be the leftmost\footnote{We consider $\gamma_{\ell+ 1}$ as a cell-path containing a single cell adjacent to the arcs $(a,b)$ and $(c,d)$ of $\partial \calQ_{k+1}$, respectively. In other words, we remove the potential section of $\gamma_{\ell+1}$ that bounces off the bottom segment. The same convention is used for $\gamma_{r-1}$ below.} cell-path in $\calQ_{k+1}$ of colour $\sigma_{\ell+1}$ from $(a,b)$ to $(c,d)$. 
	This path is revealed by using the interface exploration in the domain $\calQ_{k+1}$ with marked points $a,c,d$, started from $a$ and with $\sigma_{\ell+1}$ on its right and $1-\sigma_{\ell+1}$ on its left. More precisely, $\gamma_{\ell+1}$ is a subset of the cells on the right side of this interface. See Figure~\ref{fig:cut_tube} (left).
	We then say that step $k+1$ is successful if the endpoint of $\gamma_{\ell+ 1}$ is to the left of $z$, i.e., on the arc $(z,d)$.
	
	Similarly, if $X_{k+1} = \sfR$, let $\gamma_{r-1}$ be the rightmost path in $\calQ_{k+1}$ of colour $\sigma_{r-1}$ from $(a,b)$ to $(c,d)$.
	Then, we say that step $k+1$ is successful if the endpoint of $\gamma_{r-1}$ is to the right of $z$, i.e., on the arc $(c,z)$. 
	
	At this stage, if the exploration is successful, further reveal the state of all unexplored cells that are disconnected from the top of $\Rect(R,H)$ by the explored part. Set $\tau_{k+1}$ to be the end of this exploration step.
	If the exploration is not successful, reveal the state of all cells in arbitrary order and set $\tau_{k+1} =  \infty$. 
	
	\medskip 
	
	\noindent {\bf Induction step: checking conditions.}
	We now prove that, if the exploration is successful at step $k+1$, then the points (1), (2) and (3) are satisfied at time $\tau_{k+1}$. 
	Assume we are in the first case $X_{k+1} = \sfL$, where we explored $\gamma_{\ell+1}$ from the left. The case $X_{k+1} = \sfR$ is handled analogously, with the roles of left and right interchanged. 
    Assuming that step~$k+1$ is successful, we have that the endpoint of $\gamma_{\ell+ 1}$ is on the arc $(z,d)$.
	
	We start with point (1). 	
	Recall the tubes $T_{j}^{(k+1)}$ defined above for $ j \leq \ell$ and $j \geq r$. 
	Their exposition zones are indeed disjoint and do not intersect $\calQ_{k+1}$. As such, they are not affected by the exploration between $\tau_{k}$ and $\tau_{k+1}$, and  $e_{[\tau_{k+1}]}$ exposes the corresponding paths $\gamma_j$ for $ j \leq \ell$ and $j \geq r$.
	
	Define the tube $T_{\ell+1}^{(k+1)}$ as follows --- red tube in Figure~\ref{fig:stage_k2}.
	Follow $S_{R}$ from its start up to the first box contained in $\calQ_{k+1} \cap [-R,R] \times [H- 4 \eta_k R, H - 4 \eta_{k+1} R]$
	--- such an intersection necessarily occurs since $S_R$ exposes $\gamma_{r}$ in $e_{[\tau_k]}^c$, which is entirely below the level $H- 4\eta_k R$.
	Then use boxes contained in $\calQ_{k+1} \cap [-R,R] \times [H- 4\eta_k R, H - 4 \eta_{k+1} R]$ to connect $S_R$ to~$S_L$. 
	Finally, follow $S_L$ to its end. 
	By the fact that $\gamma_{\ell + 1}$ exits $\calQ_{k+1}$ on the arc $(z,d)$, any crossing of $T_{\ell+1}^{(k+1)}$ necessarily intersects $\gamma_{\ell+1}$.
	The construction of $T_{\ell+1}^{(k+1)}$ guarantees that the connected component of the start of $T_{\ell+1}^{(k+1)}$ in $T_{\ell+1}^{(k+1)} \setminus \gamma_{\ell + 1}$ is entirely unexplored. 
	Thus $\gamma_{\ell+1}$ is exposed by $T_{\ell+1}^{(k+1)}$ in $e_{[\tau_{k+1}]}^c$. 
	
	Furthermore, the exposition zone of $\gamma_{\ell+1}$ in $T_{\ell+1}^{(k+1)}$ is disjoint from the exposition zones of $\gamma_j$ in $T_{j}^{(k+1)}$ with $j \leq \ell$ and $j \geq r$. 
	Finally, all $T_{j}^{(k+1)}$ are of quality $\eta_{k+1}= \eta_k/2$ by definition.
	
	\begin{figure}
    	\begin{center}
    	\includegraphics[width = 0.75\textwidth, page = 2]{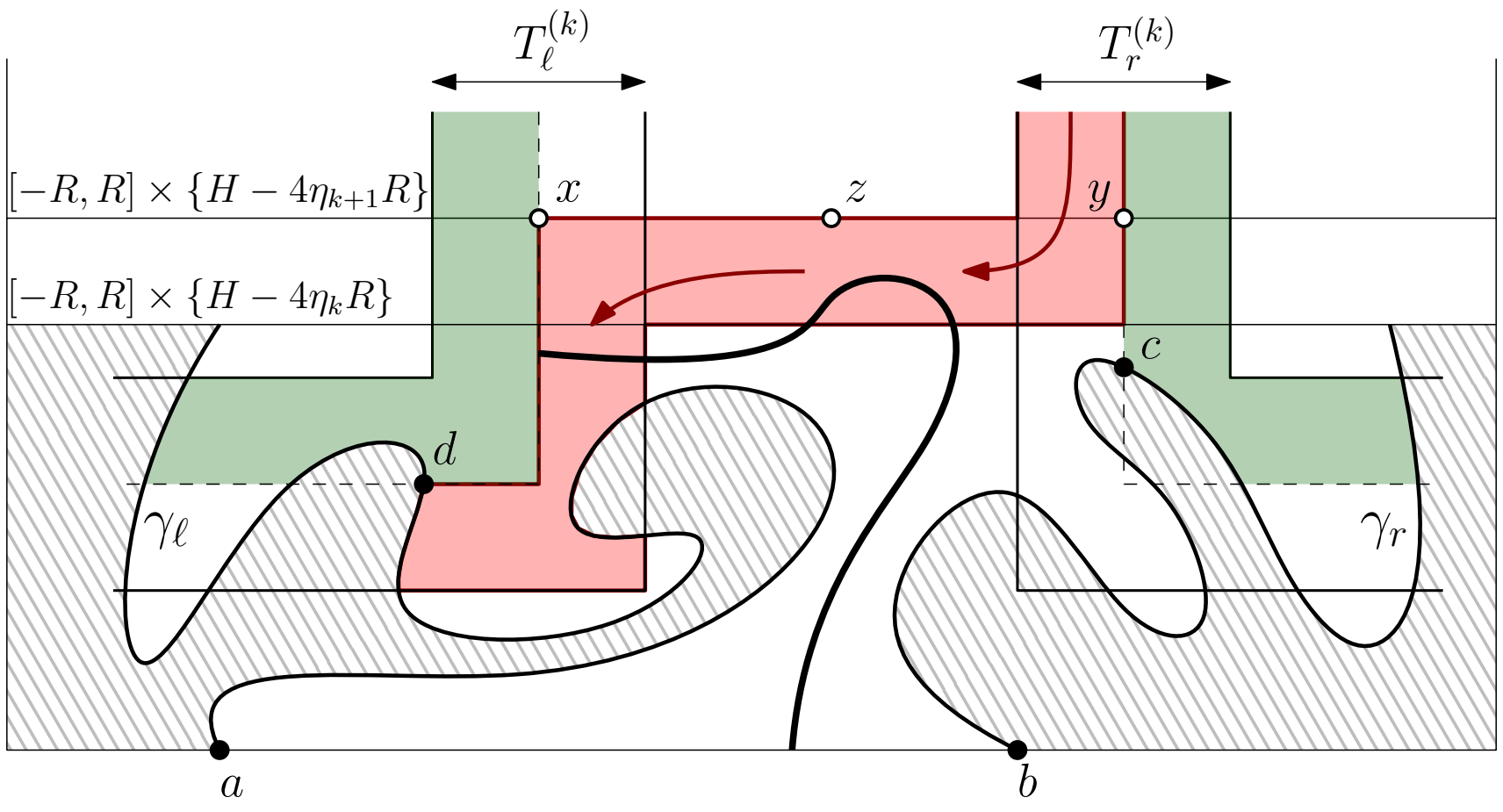}
    	\caption{A case where $X_{k+1} = \sfL$ and $\gamma_{\ell+1}$ (in thick black) lands left of~$z$. The cells explored between $\tau_k$ and $\tau_{k+1}$ are all in $\calQ_{k+1}$, left of $\gamma_{\ell+1}$. Then $\gamma_{\ell+1}$ is exposed in the unexplored part $e_{[\tau_{k+1}]}^c$ by the red tube $T_{\ell+1}^{(k+1)}$.} 
    	\label{fig:stage_k2}
    	\end{center}
	\end{figure}
	
	Point (2) is clear from the construction. Indeed, the exploration guarantees that $e_{[\tau_{k+1}]}^c$ is simply connected. 
	Furthermore, since $\calQ_{k+1}$ is contained below the line $[-R,R] \times \{H - 4\eta_{k+1} R\}$, 
	all explored cells are also contained in this part. 
	
	Finally, for point (3), we observe that, since $\gamma_{\ell+1}$ was the leftmost path of colour $\sigma_{\ell+1}$ crossing $\calQ_{k+1}$, 
	the other paths of colours $\sigma_{\ell+2}, \dots, \sigma_{r-1}$ are contained between $\gamma_{\ell+1}$ and~$\gamma_{r}$.
	
	\medskip 
	
	\noindent {\bf Conclusion:}	
    Set $\Exp = e_{[\tau_N]}$. Observe that, if $\tau_{N} < \infty$, then $\Exp$ indeed $\eta$-exposes crossings of colours $\sigma$ towards the top for $\eta = \eta_N$.
    
	\medskip 
	
	\noindent {\bf Probability of success:}	
	We turn to estimating the probability of success of the process. Observe that if $\tau_k < \infty$ then, at step $k$, for any possible configuration in the unexplored region, the leftmost crossing of $\calQ_{k+1}$ from $(a,b)$ to $(c,d)$ of colour $\sigma_{\ell+1}$ lies to the left of the 
	rightmost crossing of colour $\sigma_{r-1}$. Recall that the existence of these two crossings is a consequence of point (3).
	As such, at least one of the following holds: 
	\begin{itemize}
	\item the leftmost crossing of $\calQ_{k+1}$ from $(a,b)$ to $(c,d)$ of colour $\sigma_{\ell+1}$ exits left of $z$ or
	\item the rightmost crossing of $\calQ_{k+1}$ from $(a,b)$ to $(c,d)$ of colour $\sigma_{r-1}$ exits right of $z$.
	\end{itemize}
	We conclude that for any (undiscovered) configuration, $\tau_{k+1} < \infty$ for at least one of the cases $X_{k+1} = \sfL$ or $X_{k+1} = \sfR$. Thus
	\begin{align}
	\bbP\big[\tau_{k+1} < \infty \,|\, \tau_k < \infty, \,e_{[\tau_k]},\, \omega_{e_{[\tau_k]}}, \,A_{\sigma}({\rm Rect}(R,H))\big] \geq \tfrac12,
	\end{align}
    and by repeated applications of the above, we obtain the desired result:
    \begin{align}
        \bbP\big[\tau_N < \infty\,\big|\, A_{\sigma}({\rm Rect}(R,H))\big]\geq \big(\tfrac12\big)^N > 0.
    \end{align}
\end{proof}

\subsection{Proof of Theorem~\ref{thm:arm_exposition_rect} for general tilings}

We turn to the case of general tilings. 
Fix such a tiling $\calG$ and a percolation measure $\bbP$ on $\calG$. Moreover, fix $R,H$ and $\sigma$ as in the statement of the theorem. 
Let $\eta$ be the quantity obtained in Section~\ref{sec:favourable_tilings}.
We will prove the theorem for $\tilde\eta = \eta/4$. 

Applying the same algorithm as in Section~\ref{sec:favourable_tilings}, we obtain explorations and tubes of quality $\eta$.
Explored cells intersecting the tubes $T_1,\dots, T_N$ that are used to expose certain paths $\gamma_1,\dots, \gamma_N$ may however overlap. As such, these tubes do not satisfy the requirements of the theorem. To accommodate this, we replace each tube by a narrower tube lying strictly inside it.

The hypothesis of the theorem requires that all cells of $\calG$ intersecting $\Rect(R,H)$ have diameter at most $\tilde \eta R = \eta R/4$. 
In particular, this bound applies to every cell revealed by the exploration.
As such, each tube $T_i$ contains a tube $\tilde T_i$ of quality $\tilde \eta$ at scale $R$, 
whose start is contained in the start of $T_i$, whose end is contained in the end of $T_i$
and such that the distance between $\tilde T_i$ and $T_i^c$ is at least $\tilde \eta R$.

The algorithm above, together with the paths $\gamma_1,\dots, \gamma_N$ and the tubes $\tilde T_1,\dots, \tilde T_N$ satisfies the conclusion of the theorem. \hfill $\square$

\section{Exposition in an annulus: Proof of Theorem \ref{thm:arm_exposition}}\label{sec:proof_annulus}

Fix all quantities as in the statement of Theorem~\ref{thm:arm_exposition}.
As in the proof of Theorem~\ref{thm:arm_exposition_rect}, we will assume that our measure is already conditioned on  $A_{\sigma}(r,R)$, and thus that $\bbP[A_{\sigma}(r,R)] = 1$. 
We will see later how to choose $0 < \eta < \frac18$ and $c>0$, depending only on the length $N$ of the colour sequence. 
As stated in the theorem, we assume that all cells of $\calG$ intersecting~$\La_R$ are of diameter at most $\eta R$. 
    
Since $\sigma$ is polychromatic, there exists at least one open/closed interface crossing the annulus $\Ann(r,R)$ from $\partial\La_r$ to $\partial\La_R$. 
Our algorithm will consist of two steps.
First, we explore an interface starting from $\dt\La_r$ up to $\dt\La_{3R/4}$. This will allow us to cut the annulus in order to reduce the problem to a quad as in Theorem~\ref{thm:arm_exposition_rect}. In the second step, we apply a slight variation of the algorithm of Theorem~\ref{thm:arm_exposition_rect}. 
    
\subsection{Interface exploration in annuli}
We start by describing how to explore an interface from $\dt \La_r$ up to the first time it exits $\La_{3R/4}$. Note that, since $A_{\sigma}(r,R)$ is assumed to occur, at least $N$ cells intersect $\La_r$. After a small perturbation of $\La_r$, we may also assume that $\partial \La_r$ contains no triple point of $\calG$.
    
Fix an arbitrary ordering of the edges of $\calG$ crossing $\partial \La_r$ and orient each of them outwards. Examine these edges in the chosen order, revealing their adjacent cells if they have not been revealed previously. Do so until an edge separating opposite colours is found. Starting from this edge, follow the interface, preserving the colours on either side of the interface, until it either exits $\La_{3R/4}$ or returns to $\La_r$.
In the former case, stop the exploration. In the latter case, consider the next unexplored edge crossing $\partial \La_r$ and proceed as before. Repeat the procedure until the exploration exits $\La_{3R/4}$.
    
\begin{figure}
	\begin{center}
	\includegraphics[width = 0.5\textwidth]{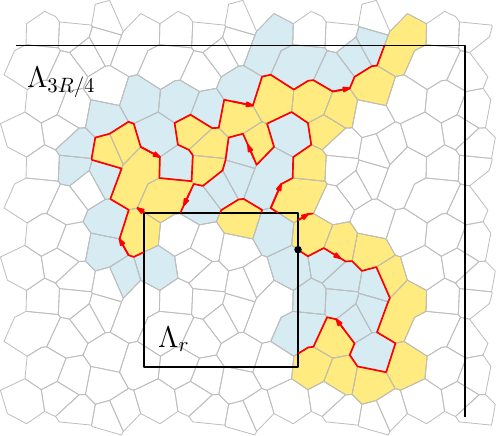}
	\caption{The explored edges, starting at the marked point on $\partial \La_r$, with the ordering being given by the counterclockwise order. The exploration process stops when reaching $\partial \La_{3R/4}$.}
	\label{fig:exploration_bouncing}
	\end{center}
\end{figure}

Due to our assumption, at least one interface exists between $\La_{r}$ and $\partial \La_R$. As such, the exploration above does eventually reach $\partial \La_{3R/4}$, and stops. We denote by $\Gamma$ the set of all edges examined during this procedure, excluding those intersecting $\partial \La_r$ which separate cells of same colour. See Figure~\ref{fig:exploration_bouncing} for an illustration. 
	
The purpose of the exploration above is to represent the arms of $A_\sigma(r,R)$ as crossings of a suitable quad to which we can adapt the argument of Theorem~\ref{thm:arm_exposition_rect}. 
The important property of the set of edges $\Gamma$ is that no arm of $A_\sigma(r,R)$ can cross $\Gamma$ since these edges are between cells of opposite colour. So one may view $\Gamma$ as a union of walls that cuts the annulus into disjoint connected components. 
There is a unique connected component of $\Ann(r,R) \setminus \Gamma$ that contains $\partial \La_R$ --- we denote this component by $\calD$.
Since the arms producing $A_\sigma(r,R)$ may not cross $\Gamma$, they all occur in $\calD$. 
Note, however, that they may use revealed cells of $\calD$ on either side of $\Gamma$. 
	
At this stage $\calD$ is topologically still an annulus, and not a rectangle as desired. 
It does however have long stretches on its inner boundary where arms may not originate from: the part of the inner boundary formed by $\Gamma$. 
We exploit this feature in the next section. 
	
\begin{rem}
	It may be tempting to actually explore the interface until it reaches $\partial \La_R$, rather than stopping at $\partial \La_{3R/4}$.
	This would produce an unexplored outer region which is a topological rectangle. 
	The problem with doing this is that the exploration may come arbitrarily close to $\partial \La_R$ several times
	before finally reaching $\partial \La_R$. 
	As such, the topological rectangle obtained above may  degenerate, 
	which would prevent the application of the algorithm of Theorem~\ref{thm:arm_exposition_rect}.
\end{rem}

\subsection{Producing a good quad}\label{sec:good_quad}

\begin{figure}
\begin{center}
	\includegraphics[height = 6cm]{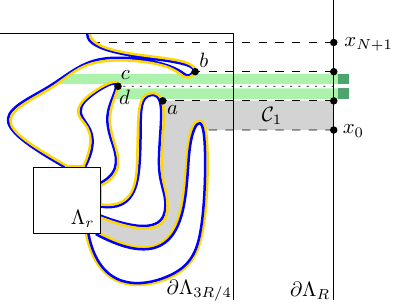}
	\hspace{ 0.02\textwidth}
	\includegraphics[height = 6cm]{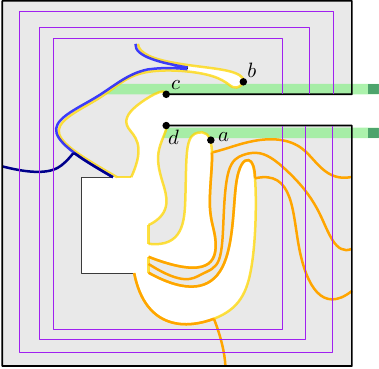}
	\caption{A sketch of the procedure in Section~\ref{sec:good_quad}.
	{\em Left:}
	The exploration of $\Gamma$ is depicted in blue and yellow. To be precise, $\Gamma$ is the edge-path between the blue and yellow sides. 
	The horizontal lines $\ell_i$ are dashed. A domain $\calC_1$ is depicted in grey. 
	We assume here that $J=2$ is good; the line~$\ell_{J-1/2}$ is dotted and the points $a,b,c,d$ are marked. 
	The exposition zones contained in $T_0^{(0)}$ and $T_{N+1}^{(0)}$ are green.
	{\em Right:} The grey region is the simply connected domain $\calD'$, with the marked points. We separated the two sides of~$\ell_{J-1/2}$ for clarity.
	The target arc $(c,d)$ (which corresponds to the top of the quad) is depicted by a slightly bolder black line. 
	Observe that the arcs $(b,c)$ and $(d,a)$ of $\partial \calD'$ correspond to the inside of $\calC_J$, which is assumed to be good, and hence no arm originates from these arcs. 
	Depicted in deeper colours are five arms, one open (dark blue) and four closed (orange). They are allowed to use the revealed cells on the boundary of $\calD'$.
	The purple lines play the role of the horizontal levels $H- 4 \eta_k R$ of Figures~\ref{fig:stage_k}~and~\ref{fig:stage_k2}.}
	\label{fig:caves}
\end{center}
\end{figure}

Recall that $\calD$ is the connected component of $\Ann(r,R) \setminus \Gamma$ containing $\dt \La_R$ and that all arms occur in this topological annulus. 
After rotating the picture, we may consider that $\Gamma$ ends on the top side of $\La_{3R/4}$. 
We now describe how to obtain a domain $\calD'$ in which we may apply the algorithm of Theorem~\ref{thm:arm_exposition_rect}. 
Figure~\ref{fig:caves} provides a useful illustration. 

Consider $N+2$ equally spaced points on the right side of $\dt \La_R$:
\begin{align}
    x_j = \bigg(R, \frac{R}{2} + \frac{jR}{4(N+1)}\bigg), \qquad j = 0, \dots, N+1.
\end{align}
Let $\ell_j$ be the line segment starting from $x_j$ and travelling horizontally to the left until it meets $\dt \calD$.
We will actually approximate $\ell_j$ by paths of edges of $\calG$; by taking $\eta$ sufficiently small, we may assume that such approximate disjoint lines exist.

For $j = 1, \dots, N+1$, let $\calC_j$ be the connected component of $\calD \setminus (\ell_{j-1}\cup \ell_j)$ that contains the part of $\dt \La_R$ between $x_{j-1}$ and $x_j$.
The regions $\calC_1,\ldots,\calC_{N+1}$ are pairwise disjoint.

We now make a crucial observation. 
Fix a configuration such that $A_\sigma(r,R)$ occurs and choose $N$ arms realising the event. Each arm has a starting point on $\dt \La_r$.
Since there are~$N+1$ regions $\calC_1, \dots, \calC_{N+1}$ and only $N$ starting points for the arms, at least one region~$\calC_j$ contains no starting point of any of the $N$ chosen arms. We say that $\calC_j$ is \emph{good} (for this configuration) 
if $A_\sigma(r,R)$  may be realised with no arms originating in $\calC_j$. 

In the final algorithm, one of the regions $\calC_j$ will be chosen uniformly at random. For now, we describe how to construct a quad in order to apply the algorithm of Theorem~\ref{thm:arm_exposition_rect} for a fixed index $J$, under the assumption that $\calC_J$ is good. 

Let $\ell_{J-1/2}$ be the horizontal line segment starting from $x_{J-1/2} = x_J - (0,\frac{R}{8(N+1)})$ and running left until the first time it meets $\dt\calD$. As before, we use a discrete approximation using edges of $\calG$.

We now set $\calD' = \calD \setminus \ell_{J-1/2}$; this is a simply connected domain. 
Define the four points (or rather prime ends to be more precise) $a,b,c,d$ in counterclockwise order on $\partial \calD'$ as follows. 
The left endpoints of $\ell_{J-1}$ and $\ell_J$ are $a$ and $b$ respectively. 
The left endpoint~$\ell_{J-1/2}$ corresponds to both $c$ and $d$; $c$ is part of the boundary as seen from above, while $d$ is seen as part of the boundary of the lower part. We denote by $(a,b)$, $(b,c)$, $(c,d)$ and $(d,a)$ the four arcs of the boundary of $\calD'$. See Figure~\ref{fig:caves}. 

Let $S = (a,b) \cap \partial \La_r$ be the subset of $(a,b)$ lying on the inner boundary $\partial \La_r$. Observe that $S$ does not interesect $\calC_J$. Thus, in any configuration for which $\calC_J$ is good, there exist arms of colours prescribed by $\sigma$ starting from $S$ and ending on the arc $(c,d)$. 
To be precise, there exists a cyclic permutation $\sigma'$ of $\sigma$ such that paths of colours $\sigma'$ cross $\calD'$ between $S$ and $(c,d)$, ordered from left to right. We denote this event by $A_{\sigma'}(\calD')$. The parts of $(a,b)$ created by the exploration of $\Gamma$ are not admissible as starting points for the arms but the arms may subsequently use cells adjacent to $\Gamma$.
The event $A_{\sigma'}(\calD')$ will serve as the starting point for the induction in the proof of Theorem~\ref{thm:arm_exposition_rect}.

Let us also define the tubes appearing in condition (1) of the induction for the initial step. 
Set $h_j = \frac{R}{2} + \tfrac{jR}{4(N+1)}$,
\begin{align}
    T_0^{(0)} &= [-(1+2\eta)R,(1+2\eta)R] \times [h_{J-1},h_{J-1/2}]     
    \qquad \text{ and } \\
    T_{N+1}^{(0)} &= [-(1+2\eta)R,(1+2\eta)R] \times [h_{J-1/2},h_J],
\end{align} 
and view them as tubes starting at their rightmost squares and ending at their leftmost squares. 
These are indeed tubes of quality $\eta_0 = \eta_0(N) = 1/(8(N+1))$. 
Observe that the tube $T_0^{(0)}$ exposes\footnote{The exposition of boundary arcs here is meant in the same sense as in footnote~\ref{footnote_curve_expo}.} the boundary arc $(da)$ of the quad $\calD'$ and the tube $T_{N+1}^{(0)}$ exposes the arc $(bc)$. 
See again Figure~\ref{fig:caves}.

Assume that we know the cyclic permutation $\sigma'$ for which $A_{\sigma'}(\calD')$ occurs. 
Then we may adapt the algorithm of Theorem~\ref{thm:arm_exposition_rect} in the domain $\calD'$ with the colour sequence $\sigma'$. Indeed, at each induction step, we select the leftmost or rightmost crossing in $\calD'$ among those starting on $S$.
The final additional detail needed is to define an equivalent for the horizontal strips 
$[-R,R] \times [H- 4 \eta_k R, H - 4 \eta_{k+1} R]$ used to switch between the left and right tube 
when constructing $T_{\ell+1}^{(k+1)}$ or $T_{r-1}^{(k+1)}$ in the induction step. 
In our case, these strips are defined as the intersections of $\calD'$ with the annuli $\Ann(R(1 - \eta_0/2^k), R(1 - \eta_0/2^{k+1}))$.

Thus, the adaptation of the algorithm, if successful, produces an exploration in $\calD'$ and disjoint tubes of quality $\eta(N) = \eta_0/2^N$ that expose paths of colours $\sigma$ towards the outside.

\subsection{Conclusion}

In conclusion, the algorithm for Theorem~\ref{thm:arm_exposition} is the following. 
\begin{enumerate}
    \item Explore $\Gamma$ starting from an arbitrary point on $\partial \La_r$ until it reaches~$\dt \La_{3R/4}$.
    \item Define the regions $\calC_j$ and choose $J$ uniformly in $\{1,\dots,N+1\}$. 
    Also, choose a cyclic permutation $\sigma'$ of $\sigma$ uniformly at random.
    \item Finally, attempt the algorithm of Theorem~\ref{thm:arm_exposition_rect} in $\calD'$ for arms of colours $\sigma'$.
\end{enumerate} 

For any fixed configuration, there exists a chance of at least $\frac{1}{N(N+1)}$ that step (2) correctly identifies a good $\calC_J$ and that $\sigma'$ is such that $A_{\sigma'}(\calD')$ occurs. 
Thus, using the estimate of Theorem~\ref{thm:arm_exposition_rect}, we obtain a uniform lower bound on the probability of success of this algorithm.
\hfill$\square$

\section{Application to FK-percolation}\label{sec:applications_FK} 

Theorems~\ref{thm:arm_exposition} and \ref{thm:arm_exposition_rect} expose arms without imposing any assumptions on the percolation measure. In this section, we show how exposition of arms can be combined with crossing estimates to obtain arm separation. We work with FK-percolation for concreteness. The argument relies on the spatial Markov property, positive association and RSW estimates (all of which will be recalled below). As explained in Section \ref{sec:comments}, the notion of exploration and the result of Theorem \ref{thm:arm_exposition} can readily be adapted to the case of bond percolation on $\bbZ^2$. We shall use such an adaptation implicitly.

\subsection{FK-percolation and the necessary tools}

Let us start by briefly recalling the definition of FK-percolation and the properties that will be used below. For more detailed statements and proofs, see \cite{Gri06, Dum20, Man25}.
Let $G = (V,E)$ be a finite subgraph of $\bbZ^2$ with vertex set $V$ and edge set $E$. A percolation configuration $\om \in \{0,1\}^E$ is the spanning subgraph of $G$ obtained by retaining the edges $e \in E$ for which $\om_e = 1$ and discarding the others. Retained edges are said to be \emph{open} while the discarded edges are said to be \emph{closed}. A boundary condition $\xi$ is a partition of the boundary vertices of $G$; vertices belonging to the same partition element are identified.

For $p \in (0,1)$ and $q > 0$, the FK-measure on $G$ with boundary condition $\xi$ is
\begin{align}
    \phi_{G,p,q}^\xi[\om] = \frac{1}{Z_{G,p,q}^\xi}p^{o(\om)}(1-p)^{|E|-o(\om)}q^{k_\xi(\om)},
\end{align}
where $o(\om)$ is the number of open edges in $\om$, $k_\xi(\om)$ is the number of connected components after the identifications prescribed by $\xi$, and $Z_{G,p,q}^\xi$ is a normalising constant. The case~$q~=~1$ corresponds to Bernoulli bond percolation.

The dual configuration $\om^*$ is defined by declaring a dual edge open precisely when the corresponding primal edge is closed. In the terminology of the previous sections, primal open paths have colour $1$ and dual open paths have colour $0$. See also the discussion in Section~\ref{sec:comments} and Figure~\ref{fig:lattices_as_tilings}.

Throughout this section, we fix $q \in [1,4]$ and $p = p_c(q) = \frac{\sqrt q}{1 + \sqrt q}$, the critical parameter~\cite{BefDum12}. We suppress these from notation. Also write $\phi$ for the unique infinite-volume FK-measure on $\bbZ^2$ \cite{DumSidTas17}. 

We will use the following form of the RSW estimates deduced readily from \cite{DumSidTas17}.  
Recall the notion of tube of quality $\eta$ at scale $R$ and that of crossing of such a tube. 

\begin{prop}
	For any $\eta >0$ there exists $c > 0$ such that, for any $R$ large enough, 
	any tube $T$ of quality $\eta$ at scale $R$ and any boundary conditions $\xi$ on $T$, 
	\begin{align}
		\tag{RSW}\label{eq:RSW}
		\begin{aligned}
			\phi_T^\xi[T \text{ contains a primal open crossing}] &\geq c \quad\text{ and } \\
		    \phi_T^\xi[T \text{ contains a dual open crossing}] &\geq c.  
     	\end{aligned}
	\end{align}
\end{prop}

Note that the crossings in the above are not required to touch the boundary of $T$; they only need to connect the starting box with the ending box of $T$. Thus, the inequalities above are indeed valid for all $q  \in [1,4]$. 

When working with \eqref{eq:RSW} estimates, a particularly convenient property of FK-percolation is \emph{positive association}: for any increasing events $A$ and $B$ and any graph $G$ and boundary conditions $\xi$
\begin{align}
    \label{eq:FKG}
    \phi_G^\xi[A \cap B] \geq \phi_G^\xi[A] \phi_G^\xi[B].
    \tag{FKG}
\end{align}
Here, an event is increasing if opening additional edges preserves its occurrence. The same inequality holds for two decreasing events. 

Finally, after an exploration, the \emph{spatial Markov property} identifies the conditional law in the unexplored region as an FK-measure with the same parameters $p,q$ and with the boundary conditions induced by the revealed configuration. 

Combining these properties allows us to construct crossings in the exposition zones and attach them to the exposed arms.

\subsection{Flower domains and separated arms}

We formulate separation of arms using flower domains, following \cite{DumMan20}. A particularly convenient feature of flower domains is that they can be explored from one side while leaving the region on the other side unexplored.

\begin{figure}
\begin{center}
\includegraphics[height = 4cm]{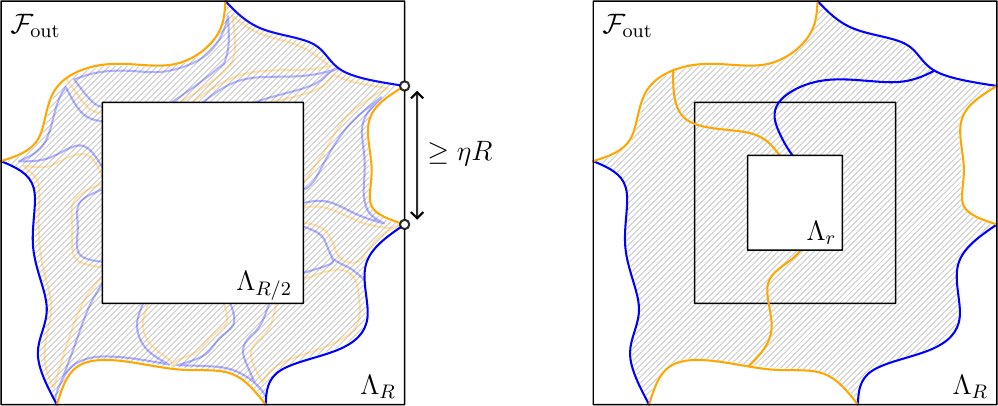}
\caption{{\em Left:} The outer flower domain between scales $\La_{R/2}$ and $\La_{R}$ is $\eta$-well-separated. This flower domain is revealed by exploring part of the hatched region. Some of the interfaces revealed during the process are drawn.
{\em Right:} A realisation of the event $A^{0,\eta}(r,R)$. This event is revealed by exploring the outer flower domain between $\La_{R/2}$ and $\La_{R}$ and the entire configuration in $\Ann(r,R/2)$. This correspond to the hatched region. In particular, the arms can be extended in the unexplored region by using the properties of FK-percolation.}
\label{fig:flower_domain}
\end{center}
\end{figure}

For $r < R$, the inner flower domain explored from $\La_R$ to $\La_r$ is obtained by exploring the interfaces of $\Ann(r,R)$ starting on $\partial \La_R$ and retaining the unexplored domain containing the inner box $\La_r$. Similarly, the outer flower domain explored from $\La_r$ to $\La_R$ is obtained by exploring the interfaces of $\Ann(r,R)$ starting on $\partial \La_r$ and retaining the domain containing the exterior of $\La_R$. See Figure~\ref{fig:flower_domain}~(left) for an illustration.

The boundaries of these domains consist of alternating primal and dual arcs, called \emph{petals}, or of a single primal or dual circuit. 
For the inner flower domain, the endpoints of petals lie on $\partial \La_r$ and for the outer flower domain, they lie on $\partial \La_R$.
An inner flower domain on $\La_r$ is \emph{$\eta$-well-separated} if it contains at least two petals and any two distinct petal endpoints are at distance greater than $\eta r$. For an outer flower domain on $\La_R$ to be $\eta$-well-separated, the corresponding requirement is that the endpoints lie at a distance greater than $\eta R$.

\begin{definition}[Separated arm events]\label{def:separation}
    Fix $\sigma \in \{0,1\}^N$ and $R \geq 4r$. Let $\Fin$ be the inner flower domain explored from $\La_{2r}$ to $\La_r$, and let $\Fout$ be the outer flower domain explored from $\La_{R/2}$ to $\La_R$.
    Let $A_\sigma^{\eta,\eta'}(r,R)$ be the event that $\Fin$ and $\Fout$ are $\eta$- and $\eta'$-well-separated, respectively, and that $A_\sigma(r,R)$ can be realised by arms connecting distinct petals of $\Fin$ to distinct petals of $\Fout$, with each arm attached to petals of its own colour.

    We write $A_\sigma^{\eta,0}$ when only the inner separation is imposed and $A_\sigma^{0,\eta'}$ when only the outer separation is imposed. See Figure \ref{fig:flower_domain} ({\em right}).
\end{definition}

Requiring distinct arms to use distinct petals also separates consecutive arms of the same colour; unused petals may lie between them.
The well-separated petals provide enough room to continue the arms through additional crossings. The following theorem shows that imposing this separation costs only a uniformly positive factor, even when the configuration away from the target boundary is fixed.

\begin{thm}[Arm separation]
    \label{thm:arm_separation}
    Fix a polychromatic sequence $\sigma \in \{0,1\}^N$. For sufficiently small $\eta, \eta' > 0$, there exists $c = c(q,\sigma,\eta,\eta') > 0$ such that the following holds for all $R \geq 4r$:
    \begin{itemize}
        \item[--] separation towards the inside:
        \begin{align}
            \phi[A_\sigma^{\eta,0}(r,R) \,|\, A_\sigma(r,R), \om_{\mathrm{out}}] &\geq c;
        \end{align}
        \item[--] separation towards the outside:
        \begin{align}
            \phi[A_\sigma^{0,\eta'}(r,R) \,|\, A_\sigma(r,R), \om_{\mathrm{in}}] &\geq c;
        \end{align}
        \item[--] separation in both directions: 
        \begin{align}
            \phi[A_\sigma^{\eta,\eta'}(r,R) \,|\, A_\sigma(r,R), \om_{\mathrm{mid}}] & \geq c.
        \end{align}
    \end{itemize}
    Here, $\om_{\mathrm{out}}, \om_{\mathrm{in}}$ and $\om_{\mathrm{mid}}$ denote fixed configurations on the respective annuli $\Ann(2r,R)$, $\Ann(r,R/2)$ and $\Ann(2r,R/2)$. Each inequality holds for every such conditioning of positive probability.
\end{thm}

To experts, the above will appear as a direct consequence of Theorem~\ref{thm:arm_exposition} and the properties of FK-percolation, particularly \eqref{eq:RSW}.
We prefer stating separation for FK-percolation in terms of flower domains, as these are explored using a standard algorithm, rather than the bespoke algorithm appearing in Theorem~\ref{thm:arm_exposition}. Furthermore, flower domains and their exploration already appeared in \cite{DumMan20, GasMan25} and seem particularly adapted to surgeries for FK-percolation and similar models.  

\begin{rem}
    In certain contexts, it may be useful to consider half-plane arm events, i.e., arm events where the arms are restricted to lie in a given half-plane. The notion of flower domain adapts naturally to this context, so one can prove an analogue of Theorem \ref{thm:arm_separation} for arm events in the half-plane by using the same method. This construction uses Theorem~\ref{thm:arm_exposition_rect} directly, rather than Theorem~\ref{thm:arm_exposition}. Thus, in the case of half-plane arm events, one is not restricted to polychromatic arm events.
\end{rem}

\subsection{From exposition to separation}
We now explain how exposed arms can be extended to well-separated arms with uniformly positive probability. Several variants of this construction will be used below, depending on the direction of exposition and the configurations already revealed.
For clarity, we state the outward version with the configuration in the inner part of the annulus revealed. 
The particular radii chosen below only serve to provide a macroscopic buffer between the revealed configuration and the regions used to construct the connections.

\begin{lem}[Extending exposed arms]
    \label{lem:extension_exposed_arms}
    Fix a sequence $\sigma \in \{0,1\}^N$ and $\eta>0$ sufficiently small. Then there exists $c = c(\sigma, \eta) > 0$ such that the following holds for all $R \geq 4r$. Fix a configuration $\om_{\mathrm{in}}$ on $\Ann(r,R/2)$ and consider the FK-measure conditioned on $\om_{\mathrm{in}}$. Suppose that an exploration $\Exp$ $\eta$-exposes paths of colours $\sigma$ towards the outside by only revealing edges in $\Ann(r,5R/8)$. Write $\om_{\mathrm{exp}}$ for the revealed configuration. Then,
    \begin{align}
        \phi[A_\sigma^{0,\eta}(r,R) \,|\, \om_{\mathrm{in}}, \om_{\mathrm{exp}}] \geq c,
    \end{align}
    as long as the conditioning is non-degenerate.
\end{lem}

\begin{proof}[Proof sketch]
    The lemma follows by a standard \eqref{eq:RSW} construction which we  sketch here. 
    Fix a realisation of the configuration $\om_{\mathrm{in}}$, and a realisation $\Exp$ of the exploration; we write~$\om_{\mathrm{exp}}$ for the explored configuration and $\Exp$ for the explored set. Write $\gamma_1,\dots, \gamma_N$ for the arms exposed by the exploration.
    
    First, we explore an $\eta$-well-separated \emph{double $2N$-petal flower domain} between $\La_{3R/4}$ and~$\La_R$. This consists of the inner flower domain $\calF_{\rm in}$ explored from $\La_{7R/8}$ to $\La_{3R/4}$ and the outer flower domain $\calF_{\rm out}$ explored from $\La_{7R/8}$ to~$\La_R$. We say that they form an $\eta$-well-separated double flower domain with $2N$ petals if 
     each has $2N$ petals, with the corresponding petals connected in $\calF_{\rm in}^c \cap \calF_{\rm out}^c$ by disjoint paths of the appropriate colours, and both are $\eta$-well-separated. See Figure~\ref{fig:double_flower_domain} for an example.
     It is a direct generalisation of \cite[Lemma 3.4]{DumMan20} that, for any $N$, there exist $\eta > 0$ and $c>0$ such that 
    \begin{align}
        \phi\big[(\calF_{\rm in}, \calF_{\rm out}) \text{ is an $\eta$-well-separated double $2N$-petal flower domain} \,\big|\, \om_{\mathrm{in}}, \om_{\mathrm{exp}}\big] \geq c,
    \end{align}
    for any configuration $\om_{\mathrm{in}}, \om_{\mathrm{exp}}$ as above.

\begin{figure}
\begin{center}
\includegraphics[width = .3\textwidth]{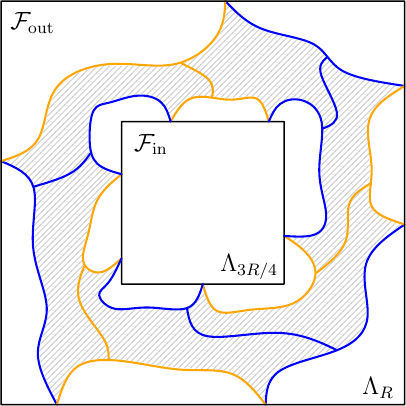}
\caption{A double $2N$-petal flower domain between $\La_{3R/4}$ and $\La_R$ for $N = 3$. This domain is revealed by exploring the corresponding inner and outer flower domains $\calF_{\rm in}$ and $\calF_{\rm out}$; the resulting explored region is hatched. Note that the petals of $\calF_{\rm out}$ and $\calF_{\rm in}$ are connected as required by the definition.}
\label{fig:double_flower_domain}
\end{center}
\end{figure}

	Assume that the exploration of $(\calF_{\rm in}, \calF_{\rm out})$ did produce an $\eta$-well-separated double $2N$-petal flower domain.
    Select $N$ distinct petals $P_1,\dots, P_N$ of $\calF_{\rm in}$ whose colours, in cyclic order, agree with $\sigma$. 
    
    By the exposition condition, there exist disjoint tubes $T_1,\dots, T_N$ of quality $\eta'$ at scale~$R$ 
    such that, for each $k$, any crossing of $T_k$ induces a connection in $\calF_{\rm in} \cap \Exp^c$ between $\gamma_k$ and $P_k$. 
    Here $\eta' > 0$ is a constant depending only on $\eta$ and $N$.
    
   Assume now that $k$ is such that $\sigma_k = 1$. 
   Conditionally on $\om_{\mathrm{in}}, \om_{\mathrm{exp}}$, the exploration producing $\calF_{\rm in}$, 
   and any configuration outside of $T_k$, 
   the measure inside $T_k \cap \calF_{\rm in} \cap \Exp^c$ dominates $\phi_{T_k}^0$. 
   Indeed, the only edges on the boundary of this domain that are contained in $T_k$ are those of $\gamma_k$ or $P_k$, which are open. 
   
   As such, we find 
   \begin{align}
        \phi\big[\,\text{each $\gamma_k$ connects to $P_k$ by a $\sigma_k$-type path  in $\calF_{\rm in} \cap \Exp^c$} \,\big|\,  \om_{\mathrm{in}}, \om_{\mathrm{exp}}, \calF_{\rm in}, \calF_{\rm out} \big] \geq c_0^N,
    \end{align}
    with $c_0$ given by \eqref{eq:RSW}. 
    Combining the two inequalities above yields
	\begin{align}
        \phi[A_\sigma^{0,\eta}(r,R) \,|\, \om_{\mathrm{in}}, \om_{\mathrm{exp}}] \geq c.c_0^N,
    \end{align}
    which concludes the proof.
   \end{proof}

\begin{proof}[Proof of Theorem~\ref{thm:arm_separation}]
    We first prove separation towards the outside; separation towards the inside can be treated similarly. Fix a configuration $\om_{\mathrm{in}}$ in $\Ann(r,R/2)$ and assume that $\phi[A_\sigma(r,R) \,|\, \om_{\mathrm{in}}] > 0$. Apply Theorem~\ref{thm:arm_exposition} in the annulus $\Ann(r,5R/8)$ to the arm event corresponding to $\sigma$. This yields an exploration $\Exp$ that reveals no edges outside $\La_{5R/8}$ and, on success, exposes paths towards the outside of $\Ann(r,5R/8)$ with some quality $\eta > 0$. Write $\om_{\mathrm{exp}}$ for the configuration revealed by the algorithm. Since the exposition occurs with positive probability uniformly over the conditioning on~$\om_{\mathrm{in}}$, we have
    \begin{align}
        \phi[\Exp \text{ $\eta$-exposes paths of colours $\sigma$} \,|\, A_\sigma(r,R), \om_{\mathrm{in}}] \geq c_{\mathrm{exp}},
    \end{align}
    with $c_{\mathrm{exp}} = c_{\mathrm{exp}}(\eta,N) > 0$ being the same constant as in Theorem~\ref{thm:arm_exposition}. 

    Let us write $\calE$ for the event that $\Exp \text{ $\eta$-exposes paths of colours $\sigma$}$. We now work under the measure $\phi[\, \cdot \, | \, \calE, \om_{\mathrm{in}}]$, without conditioning on the arm event. Lemma~\ref{lem:extension_exposed_arms} gives
    \begin{align}
        \phi[A_\sigma^{0,\eta}(r,R)\,|\, \calE, \om_{\mathrm{exp}}, \om_{\mathrm{in}}] \geq c_{\mathrm{ext}}
    \end{align}
    for any possible realisation of the explored configuration $\om_{\mathrm{exp}}$.
    Averaging over successful explorations, it follows from the previous inequalities that
    \begin{align}
        \phi[A_\sigma^{0,\eta}(r,R) \,|\, \om_{\mathrm{in}}] \geq c_{\mathrm{ext}}\, \phi[\calE \,|\, \om_{\mathrm{in}}] \geq c_{\mathrm{ext}}\,c_{\mathrm{exp}} \, \phi[A_\sigma(r,R)\,|\, \om_{\mathrm{in}}].
    \end{align}

    Internal separation follows by the same argument using exposition and extension towards the inside instead.

    \medskip

    We finally treat the case of separation in both directions. There is a slight subtlety here in that the outward exploration must preserve the distinct attachments to the inner flower domain obtained by the inward separation or vice versa. The following argument shows how the exploration algorithm can be adapted to deal with this issue. 
    
    Fix a configuration $\om_{\mathrm{mid}}$ in $\Ann(2r,R/2)$. Averaging the separation estimate towards the inside over configurations outside $\La_{2r}$ gives
    \begin{align}
        \phi[A_\sigma^{\eta,0}(r,R) \,|\, \om_{\mathrm{mid}}] \geq c_{\mathrm{in}}\, \phi[A_\sigma(r,R) \,|\, \om_{\mathrm{mid}}].
    \end{align}
    
    Fix now a configuration $\tilde \om$ in $\Ann(r,R/2)$ such that $\phi[\tilde \om \,|\, A_\sigma^{\eta,0}(r,R), \om_{\rm mid}] > 0$.
    We now wish to apply the algorithm of Theorem \ref{thm:arm_exposition}  for the measure $\phi[\,\cdot\,|\, A_\sigma^{\eta,0}(r,R), \tilde\om]$. 
    While this measure is conditioned on having arms separated towards the inside, it may still be possible to find a family of arms that are not separated towards the inside. Thus, we have no guarantee that the algorithm of Theorem \ref{thm:arm_exposition} exposes toward the outside a familly of arms that are already separated toward the inside. 
    
   	To overcome this small difficulty, we adapt the algorithm of Theorem~\ref{thm:arm_exposition} to only search for arms originating from distinct petals of $\calF_{\rm in}$. 
	More precisely, change the definitions of $\gamma_{\ell+1}$ and $\gamma_{r-1}$ by taking $\gamma_{\ell+1}$ (resp.\ $\gamma_{r-1}$) to be the leftmost (resp. rightmost) path crossing the quad $\calQ_{k+1}$ among those starting from an unused petal of $\Fin$.

    Using this adapted algorithm and the same argument as for separation towards the outside gives
    \begin{align}
        \phi[A_\sigma^{\eta,\eta'}(r,R)\,|\,  A_\sigma^{\eta,0}(r,R), \tilde\om, \om_{\rm mid}] \geq c_{\mathrm{out}}.
    \end{align}
    Averaging over $\tilde \om$ and combining the two inequalities gives
    \begin{align}
        \phi[A_\sigma^{\eta,\eta'}(r,R) \,|\, \om_{\mathrm{mid}}] \geq c_{\mathrm{out}} \, c_{\mathrm{in}}\, \phi[A_\sigma(r,R)\,|\, \om_{\mathrm{mid}}],
    \end{align}
    which concludes the proof.
\end{proof}

\subsection*{Statement on the use of AI} 

This work was entirely developed without the use of AI.
AI was then used in proofreading selected passages, which were then rewritten manually.

\subsection*{Acknowledgements} This work was financially supported by the Swiss National Science Foundation (SNSF). The third author further acknowledges support from the swissuniversities ``cotutelle de th\`{e}se'' grant. 

\bibliographystyle{alpha}
\bibliography{biblio.bbl}
    
\end{document}